\documentclass[a4paper,11pt]{amsart}
\usepackage{amssymb}
\usepackage[dvips]{graphicx}
\usepackage{amscd}
\newcommand{\comm}[1]{}

\def\dist{\operatorname{dist}}

\def\ti{\tilde}

\def\({\left(}
\def\){\right)}
\def\oli{\overline}

\def\raw{\rightarrow}

\def\no={\neq}
\def\sm{\setminus}

\def\C{{\mathbb C}}
\def\D{{\mathbb D}}

\def\N{{\mathbb N}}

\def\P{{\mathbb P}}
\def\Q{{\mathbb Q}}

\def\CC{{\mathcal C}}
\def\DD{{\mathcal D}}

\def\FF{{\mathcal F}}

\def\JJ{{\mathcal J}}

\def\MM{{\mathcal M}}
\def\NN{{\mathcal N}}
\def\OO{{\mathcal O}}

\def\SS{{\mathcal S}}

\def\al{\alpha}
\def\be{\beta}
\def\ga{\gamma}
\def\de{\delta}

\def\vep{\varepsilon}

\def\th{\theta}

\def\la{\lambda}

\def\om{\omega}

\def\La{\Lambda}

\theoremstyle{plain}
\newtheorem{Main}{Theorem}

\newtheorem{Thm}{Theorem}[section]
\newtheorem{Prop}[Thm]{Proposition}
\newtheorem{Lem}[Thm]{Lemma}

\newtheorem*{Claim}{Claim}

\theoremstyle{remark}
\newtheorem{Rem}[Thm]{Remark}
\newtheorem{Def}[Thm]{Definition}

\begin{document}
\begin{center}
\end{center}
\title[Rational Misiurewicz maps are rare II]
{Rational Misiurewicz maps are rare II}
\author{Magnus Aspenberg}
\address{Universit\'e Paris-Sud, Laboratoire de Math\'ematique, B\^atiment 425, UMR 8628, 91405 Orsay, France}
\email{Magnus.Aspenberg@math.u-psud.fr}
\thanks{The author gratefully acknowledges funding from the Swedish Research Council, for his stay at Laboratoire de Math\'ematique at Universit\'e Paris-Sud, Orsay. }
\begin{abstract}
We show that the set of Misiurewicz maps has Lebesgue measure zero in the space of rational functions for any fixed degree $d \geq 2$.
\end{abstract}

\maketitle


\section*{Introduction}
The notion of Misiurewicz maps has its origin from the paper \cite{MM}
from 1981 by M. Misiurewicz. The (real) maps studied in this paper have, among other things, no sinks and the omega limit set $\om(c)$ of every critical point $c$ does not contain any critical point.
In particular, in the quadratic family $f_a(x)=1-ax^2$, where $a \in
(0,2)$, a Misiurewicz map is a non-hyperbolic map where the critical
point $0$ is non-recurrent. D. Sands showed in 1998 \cite{DS} that these maps has Lebesgue measure zero, answering a question posed by Misiurewicz in \cite{MM}.
In this paper we state a corresponding theorem for rational maps on the Riemann sphere.
On the contrary, the set of Misiurewicz maps has full Hausdorff dimension (i.e. equal to the dimension of the parameter space), see \cite{AG}.

In the complex case, there have been some variations on the definition
of Misiurewicz maps. (The sometimes used definition of being a
postcritically finite map is far too narrow to adopt here.)
In \cite{SvS}, S. van Strien studies Misiurewicz maps with a
definition similar to the definition in \cite{MM}, (allowing
super-attracting cycles but no sinks). In \cite{GKS} by J. Graczyk,
G. \'Swiatek, and J. Kotus, a Misiurewicz map is roughly a map for
which every critical point $c$ has the property that $\om(c)$ does not
contain any critical point, (allowing sinks but not super-attracting
cycles). In this paper we allow attracting cycles, and only care about critical points on the Julia
set (suggested by J. Graczyk).

Let $f(z)$ be a rational function of a given degree $d \geq 2$ on the
Riemann sphere $\hat{\C}$. Let $Crit(f)$ be the set of critical points
of $f$, $J(f)$ the Julia set of $f$ and $F(f)$ the Fatou set of
$f$.

\begin{Def}
A {\em Misiurewicz} map $f$ is a non-hyperbolic rational map that has no parabolic cycles and such that $\om(c) \cap Crit(f) = \emptyset$ for every $c \in Crit(f) \cap J(f)$.
\end{Def}

We prove the following.
\begin{Main}
The set of Misiurewicz maps has Lebesgue measure zero in the space of rational functions for any fixed degree $d \geq 2$. \label{miserere}
\end{Main}

The Misiurewicz maps are a special type of Collet-Eckmann
maps, which - on the contrary - have positive Lebesgue measure in the parameter space of rational maps of any fixed degree $d \geq 2$, (see
\cite{MA}). Hence, Theorem \ref{miserere} shows that for a typical
(non-hyperbolic) Collet-Eckmann map, the critical set is recurrent.

\subsection{Two definitions of Misiurewicz maps.}
Let $\MM^d$ be the set of Misiurewicz maps of degree $d$ according to the definition above. Define
\[
P^k(f,c) =  \overline{ \bigcup_{n > k} f^n(c)} \text{ and } P^k(f) =  \bigcup_{c \in Crit(f)} P^k(f,c)  .
\]
Set $P^0(f,c)=P(f,c)$. The set $P(f)$ is the {\em postcritical set} of
$f$. We will also use the notion postcritical set for $P^k(f)$ for some suitable $k \geq 0$. Let $SupCrit(f)$ be the set of critical points in super-attracting cycles.
For $B(z,r)= \{ w :  |w-z| < r \}$, let
\[
U_{\de} = \bigcup_{c \in Crit(f) \sm SupCrit(f)} B(c,\de).
\]


We introduce another equivalent definition of Misiurewicz maps as
follows.
\begin{Def}
Let $\ti{\MM}^d$ be the set of non-hyperbolic rational maps of degree
$d$ without parabolic periodic points such that for which every critical
point $c$ satisfies $\om(c) \cap (Crit(f) \sm SupCrit(f)) =
\emptyset$. Let
\begin{align}
M_{\de,k} = \{ f \in \ti{\MM}^d: P^k(f) \cap U_{\de} = \emptyset \}
\label{dek-misse}
\text{ and } M_{\de} = \cup_{k \geq 0} M_{\de,k}.
\end{align}
\end{Def}
Then it is easy to see that $\ti{\MM}^d=\MM^d$.
If $f \in M_{\de,k}$ then we say that $f$ is {\em $(\de,k)$-Misiurewicz}. If $f \in M_{\de}$ then we say that $f$ is {\em $\de$-Misiurewicz }.
So for every Misiurewicz map $f$, there are constants $\de > 0$ and $k \geq 0$ such that $P^k(f) \cap U_{\de} = \emptyset$.


\subsection{Good families of rational maps.}
To show Theorem \ref{miserere} we will divide the parameter space of
rational functions into so called ``good'' analytic families. These
families are in a sense rigid meaning that they cannot be
quasiconformal conjugacy classes of Misiurewicz maps, unless they are
conformal conjugacy classes.
The idea consists firstly of fixing the multipliers of every present attracting cycle and was suggested by J. Graczyk. The construction of good families has also a strong resemblance with
the work by Ma\~n\'e, Sad, Sullivan \cite{MSS}, Theorem E and the proof thereof.
In every such good family, we prove a one-dimensional slice-version of
Theorem \ref{miserere}, namely Theorem \ref{discthm} (see below).

Every rational map of degree $d$ can be written in the form
\begin{equation}
R(z)=\frac{P(z)}{Q(z)}=\frac{a_0+a_1z+\ldots+a_dz^d}{b_0+b_1z+\ldots+b_dz^d},
\label{rat}
\end{equation}
where $a_d$ and $b_d$ are not both zero. Without loss of generality we may assume that $b_d = 1$. The case $a_d \neq 0, b_d = 0$ is treated analogously. Hence, the set of rational functions of degree $d$ is a $2d+1$-dimensional complex manifold and subset of the projective space
$\P^{2d+1}(\C)$. Now, simply take the measure on the coefficient space
in one of the two charts $a_d=1$ or $b_d=1$. There also is a
coordinate independent measure on the space of rational maps of a
given degree $d$, induced by the Fubini-Study
metric (see \cite{FS}). The Lebesgue measure on any of the two charts
is mutually absolutely continuous to the Fubini-Study measure.

A family of rational maps $R_a$ for $a \in V \subset \C^m$, where $V$ is open
and connected, is {\em normalized} if any two functions $R_a$ and $R_b$, $a,b
\in V$, are conformally conjugate if and only if $a = b$. If $f$ and $g$ are conformally conjugate then they are conjugate by a M\"obius transformation
\[
T(z) = \frac{\al + \be z}{\ga + \de z}.
\]
The set of M\"obius transformations forms a $3$-dimensional complex
manifold. Introduce an equivalence relation $\sim$ on the parameter
space, saying that $f \sim g$ if and only if $f = T^{-1} \circ g \circ
T$, for some M\"obius transformation $T$. Every equivalence class is a
complex $3$-dimensional manifold. These manifolds form a foliation of
the space of rational functions of degree $d$ (see e.g. Frobenius
Integrability Theorem in \cite{MiSp}). Hence to prove Theorem
\ref{miserere}, by Fubini's Theorem, it suffices to consider families
of normalised maps. If fact, we will consider $1$-dimensional slices
of normalised so
called {\em good} families in a neighbourhood of a starting function $R=R_0$, where we
fix the multiplier of every present attracting cycle and a little
more, as follows. The reader may also look at \cite{MSS}, the proof of Theorem E, for a similar construction.

First, we want to avoid the situation when critical points split under
perturbation. This is however a rare event in the parameter
space. Indeed, by classical theory (see e.g. \cite{Chirka}, Theorem p. 7) the set where critical points of higher multiplicity occurs is an analytic (discriminant) set, which has codimension $1$.
Hence we can assume that all critical points are non-degenerate, i.e. they do not split.

Now we turn to the construction of good families. For every attracting periodic point $x(a)$ of period $p$ we put
\[
F(a) = (R^p)'(x(a),a)-\la,
\]
where $\la = (R^p)'(x(0),0)$ is the multiplier of the
attracting cycle for the starting function $R=R_{0}$ and $R_{a}(z)=R(z,a)$, $a \in
\C^{2d-2}$. Then by classical theory (see e.g. \cite{Chirka}) the set $F^{-1}(0)$ is an analytic set of codimension $1$. Moreover, (see e.g. \cite{JM}, Lemma 1, p. 11) we have that $F^{-1}(0)$ is a submanifold $DF(0) \neq 0$.
Hence the set where we have fixed all the multipliers of every present attracting cycle is an analytic set of codimension equal to the number of attracting cycles.

After fixing the multipliers we get an analytic set $M$ where any two functions $f,g \in M$ are conformally conjugate in a neighbourhood of the attracting cycles. Now we want to possibly reduce $M$ further, so as to obtain a new analytic set $N \subset M$ on which this local conformal conjugacy can be extended to the postcritical set in the Fatou set. The linearization function (in the geometrically attracting case) is unique up to a multiplication by a constant and the B\"ottcher function (in the super-attracting case) is unique up to a multiplication by an $n$th root of unity. Therefore we do not need to reduce the manifold further if there is only one critical point in every basin (including the critical point sitting in the periodic orbit in the super-attracting case). Indeed, assume that $V$ is a neighbourhood of a geometrically attracting fixed point where $R_{a}$ is linearizable. Assume that we have fixed the multiplier of this fixed point. There exists a conformal map $\varphi_{a}$, analytically depending on the parameter $a$, such that
\[
\varphi_{a} \circ R_{a} (z) = \la \varphi_{a} (z),
\]
for $z \in V$, where $|\la| < 1$. Let $c(a)$ be a critical point in the basin. Hence for some $n$, $R_{a}^n(c(a)) \in V$. Since $\varphi_{a}$ is unique up to a multiplicative constant, we can adjust $\varphi_{a}$, (replacing it by $C(a) \varphi_{a}$ for some analytic function $C(a)$), so that $\varphi_a(R_a^n(c(a)))$ is constant in a neighbourhood of $a=0$. However, to fix $\varphi_a(R_a^n(c_1(a)))$ for some second critical point $c_1(a)$ we get a new analytic set and the dimension drops by one. To see this, consider the map
\[
\psi_n(a)= \varphi_a(R_a^n(c_1(a))): \C^{2d-2} \raw \C.
\]
Let $\sim_{\la}$ be the equivalence relation in $\C$ given by the action of the map $z \raw \la z$, i.e. $z \sim_{\la} \la z$, where $\la$ is the multiplier of the attracting cycle. Put $M_1 = M/ \sim_{\la}$. Then $M_1$ is homeomorphic to a torus and $\psi_{n+1}(a)=\la \psi_{n}(a)$ so we can define $\psi=\psi_n(a)$ on $M_1$, by $\psi(a)=\psi_n(a)$ for some $n$ satisfying $R_a^n(c_1(a)) \in V$. The set $\psi^{-1}(z_0)$, for some $z_0 \in M_1$ where $z_0 \sim_{\la} \varphi_0(R_{0}^n(c_1(0)))$ is a new analytic set of codimension $1$.

Hence, if there are more than one critical point in a basin of an attracting
cycle then for every surplus critical point we choose an analytic set where the (possibly iterated) critical value of this critical point
stays constant inside the conjugating domain. We are left with
an analytic set $N \subset M$, where every local conformal conjugacy
around the attracting cycles can be extended to the postcritical set. It is now straightforward to extend this local conjugacy to whole basin of attraction of every attracting cycle by a standard pullback argument. Indeed, assume that we have a conjugacy $\varphi$ on $U$, where $U$ is the union of a neighbourhood of the postcritical set in the Fatou set for $R=R_0$ and a neighbourhood around each attracing cycle:
\[
\varphi \circ R_0(z) = R_a \circ \varphi(z), \text{ for all $z \in U$}.
\]
So let us define $\varphi$ on $R^{-1}(z)$, where $z \in U$ and $R^{-1}(z) \notin U$. Since $z$ is not a critical value there are precisely $d$ different preimages $z_1, \ldots, z_d \in R^{-1}(z)$. We shall define $\varphi$ on each $z_i \notin U$. The functional equation above gives $\varphi(z_i) = R_a^{-1} \circ \varphi \circ (z)$, since $z=R_0(z_i)$. The fact that $z_i \notin U$ means that if the perturbation $r > 0$ is sufficiently small there is one point $x_i \in R_a^{-1} \circ \varphi(z)$ which is closest to $z_i$. Define $\varphi(z_i)=x_i$. Continuing in this way defines $\varphi$ in the full basin of attraction.

Let us calculate the dimension of $N$. Every geometrically attracting cycle must have at least one critical point in its basin. Therefore the total number of equations of fixing multipliers and fixing iterated critical values coincides with the number of critical points in the Fatou set. Apparently, the
dimension of $N$ is the dimension of the full parameter space, i.e. $2d-2$, reduced by the number of the critical points in the Fatou set. Hence $dim(N)$ is equal to the number of critical points in the Julia set for $R$. Let us call the connected component of such an analytic set containing $R$, a {\em good family} around $R$. Note that a good family is not necessarily unique (one can take any analytic subset of a good family, and this will form a new good family). However, we can speak of a unique {\em maximal} good family for $R$, if it is not contained in any other non-trivial good family for $R$. We now demonstrate an important property of good families.
\begin{Lem}
For any two good families $N_1$ and $N_2$ we must have that either $N_1 \cap
N_2 = \emptyset$, $N_1 \subset N_2$ or $N_2 \subset N_1$.
\end{Lem}
\begin{proof}
Recall that any good family is associated to a rational map $f$. Let
$f_i$ correspond to $N_i$, $i=1,2$. Every attracting cycle or critical
point in the Fatou set for $f_1$
and $f_2$ gives rise to analytic sets being zero sets
of a finite number of equations:
\begin{align}
N_1 &= \{a : e_1(a)=0, e_2(a)=0, \ldots, e_{k_1}(a)=0 \}, \nonumber \\
N_2 &= \{a : f_1(a)=0, f_2(a)=0, \ldots, f_{k_2}(a)=0 \}, \nonumber
\end{align}
where each $e_j(a)=0$ or $f_j(a)=0$ is an equation corresponding to either fixing a
multiplier of an attracting cycle or fixing an iterated critical value
for a critical point in the Fatou set. It is clear that if there is a
solution to $e_i(a)=f_j(a)$ for some $i$ and $j$ then $e_i \equiv
f_j$. Hence $\{f_1,\ldots, f_{k_2}\} \subset
\{e_1,\ldots,e_{k_1} \}$ if and only if $N_1 \subset N_2$. Finally, if
there are some $e_i$ and $f_j$ such that
\[
e_i,f_j \notin \{ e_1,\ldots, e_{k_1} \} \cap \{ f_1, \ldots, f_{k_2} \}
\]
then we must have that any $a \in N_1$ is thus a solution to
$e_i(a)=0$ but not $f_j(a)=0$. Hence $N_1 \cap N_2 = \emptyset$. This
completes the proof of the lemma.
\end{proof}

We will consider $1$-dimensional slices of any good family $N$. Moreover, since analytic sets are complex manifolds almost everywhere, we can assume that $N$ is a complex manifold. The following theorem is the main object of this paper.
\begin{Main} \label{discthm}
Assume that $R_a$, $a \in \C$, is an analytic normalized good family of
rational maps in a neighborhood of $a=0$ and that $R_0$ is a Misiurewicz map.
Then the Lebesgue density at $a=0$ is strictly less than $1$ in the
set of $(\de,k)$-Misiurewicz maps for any $\de > 0$ and $k \geq 0$.
\end{Main}

\begin{Rem}
Note that for example that Theorem \ref{discthm} applies to the the family $f_c(z)=z^d+c$, for any $d \geq 2$.
\end{Rem}

\begin{proof} [Proof that Theorem \ref{discthm} implies Theorem \ref{miserere}.]
We begin to prove that the set of $(\de,k)$-Misiurewicz maps have measure zero in any good family $N$ for any $\de > 0$ and $k \geq 0$. Let $E$ be the set of points in $N$ for which $N$ is a good family.

Theorem \ref{discthm} implies that no $(\de,k)$-Misiureiwcz point in any good one-dimensional slice of $N$ has Lebesgue density $1$. Make a foliation of the set $N$ into one-dimensional slices $\SS$. Since $N$ is a good family for every map in $E \subset N$, every such one-dimensional slice $S \in \SS$ is a good $1$-dimensional slice for any point $x \in E \cap S$. Now apply Theorem \ref{discthm} to each $S \in \SS$. By Fubini's Theorem we get that the set of $(\de,k)$-Misiurewicz maps in $E$ has $dim(N)$-dimensional Lebesgue measure zero.

The other good families inside $N$ which has lower dimension than $dim(N)$ are treated in a similar way. Suppose that $N$ is determined by the equations $\{e_1(a)=0,\ldots,e_n(a)=0\}$. Fix some $d < dim(N)$ and take some good subfamily $N' \subset N$ which has $dim(N')=d$. Now, since $N' \subset N$ and the inclusion is proper, $N'$ is determined by a solution set of a finite list of equations, $\{e_1(a)=0, \ldots, e_{n'}(a)=0 \}$, where $n' > n$. Put $k=n'-n$. The sets $$N(w)=\{e_1(a)=0,\ldots,e_n(a)=0,e_{n+1}(a)=w_1,\ldots,e_{n'}(a)=w_k\}$$ all belong to $N$ and are analytic sets for each $w=(w_1,\ldots,w_k) \in \C^k$
in a neighbourhood of $(0,\ldots,0) \in \C^k$.
This means that locally around almost any given point $x \in N'$ there is a foliation of $N \cap U$, where $U$ is a neighbourhood of $x$ in $N$ and where each leaf is a set of the form $N(w)$ (apart from the set of singularities which is an analytic set of measure zero in each good familiy).

Take some leaf $N'$ of this foliation and let $E' \subset N'$ be the set of $(\de,k)$-Misiurewicz maps for which $N'$ is a good family. By the same argument as above (Theorem \ref{discthm} and Fubini's Theorem), the set $E'$ has $dim(N')$-dimensional Lebesgue measure zero.
Because of the foliation structure of the sets $N' \subset N$, by Fubini's Theorem it follows that the $dim(N)$-dimensional Lebesgue measure is zero for all $(\de,k)$-Misiurewicz maps in $N$ whose good family has dimension $d$.

Continuing in this way for all dimensions $d < dim(N)$, we get finally that the set of $(\de,k)$-Misiurewicz maps has $dim(N)$-dimensional Lebesgue measure zero in $N$.

Now note that the full parameter space is decomposed into maximal good
families, which again may have different dimension.
These maximal families again form a foliation of the parameter space, locally. Applying Fubini's Theorem in some neighbourhood of any point
in a maximal family in the same way as above, we get that the full $(2d-2)$-dimensional Lebesgue measure of the set of $(\de,k)$-Misiurewicz maps is equal to zero. We arrive at Theorem \ref{miserere} noting that
\[
\mu \biggl( \bigcup_{n,k \in \N} M_{1/n,k} \biggr) \leq \sum_{n,k \in \N} \mu(M_{1/n,k}) = 0,
\]
where $\mu$ is the Lebesgue measure.
\end{proof}

There are some similarities between the methods in this paper and the paper \cite{DS} by D. Sands. The existence of a continuation
of the postcritical set in the real case in \cite{DS} is replaced by a
similar idea, namely that of a holomorphically moving postcritical set
in the complex case. Similar ideas appears in \cite{SvS}, by S. van Strien. This paper uses much of the ideas in \cite{MA} and some fundamental results from the paper by M. Benedicks and L. Carleson \cite{BC2} (and \cite{BC1}).

\subsection*{Acknowledgements}
I am grateful to Michael Benedicks for many valuable comments and discussions, especially on the transversality condition.
I am grateful to Jacek Graczyk for communicating the ideas on a suggested significant generalisation of the original result, where no sinks were allowed. I wish to thank Duncan Sands, Nicolae Mihalache and Neil Dobbs for interesting remarks and discussions on a preliminary version.
I am thankful for very interesting discussions with Dierk Schleicher at an early stage of this paper. I want to express my gratitude to the referee for many good and useful remarks and suggestions of improvements.

Finally, I want to express my warm thanks to Nan-Kuo Ho for useful comments and her encouragement during the writing of this paper. I dedicate this paper to her.
I gratefully acknowledge the hospitality of the Department of Mathematics at
Universit\'e Paris-Sud, where this paper was written.
\section{Some definitions and proof outline}\label{outlines}




Consider a one dimensional good normalised slice and take a ball $B(0,r)$ in this slice of radius $r > 0$. Let $R_0(z)=R(z)=P(z)/Q(z)$ be the starting unperturbed rational map of degree $d=\max{(\deg(P),\deg(Q))}$ and assume that $R_0(z)$ is Misiurewicz. We will study a critical point $c=c(a)$ dependent on the parameter $a \in B(0,r)$, for some (sufficiently small) $r > 0$. We sometimes write $R(z,a) = R_a(z)$. Put
\[
\xi_{n,j}(a)=R^n(v_j(a),a),
\]
where $c_j(a) \in Crit(R_a)$, $v_j(a)=R^{k_j}(c_j(a),a)$ and where $k_j$ is chosen so that
$v_j(0)=R^{k_j}(c_j(0),0)$ has no critical points in its forward
orbit. (A priori there can be finite chains of critical
points mapped onto each other. Therefore we assume that $v_j(a)$ is
the last critical value).
For simpler notation, we sometimes drop the index $j$ and write only $\xi_{n,j}(a) = \xi_n(a)$.

We also make the following convention. Chosen $\de > 0$, we always
assume that the parameter disk $B(0,r)$ is chosen so that the critical
points $c_i(a)$ move inside $B(c_i,\de^{10})$ as $a \in B(0,r)$.


We use the spherical metric and the spherical derivative unless otherwise stated.

\comm{
\begin{Def}
Given two complex numbers $A$ and $B$, we write $A \sim B$
meaning that there is a constant $C > 0$ depending only possibly on the unperturbed function $R$, $\de'$, and the perturbation
$r$ and such that the following holds:
\[
\biggl| \frac{A}{B} - 1 \biggr| \leq C.
\]
Moreover, we require that for any $\vep > 0$ there exist $\de',r > 0$ such that $C \leq \vep$.

If $A$ and $B$ are real and positive and $A \geq B - C$, then we write $A \gtrsim B$, if in addition for any $\vep > 0$ there exist $\de',r > 0$ such that $C \leq \vep A$. In particular $B \gtrsim A$ and $A \gtrsim B$ implies $A \sim B$.

\end{Def}
}


The proof consists of taking an arbitrarily small parameter ball $B(0,r)$ where $\xi_n(B(0,r))$ grows
to the large scale $S$, such that we have control of the shape of
$\xi_n(B(0,r))$. However,
$\xi_n(B(0,r))$ is not necessarily injective and we cannot expect to
have bounded distortion in the whole ball $B(0,r)$. Instead we show,
using strong argument distortion exstimates developed in \cite{MA},
that $\xi_n(B(0,r))$ ``grows nicely'' up to the large scale so that
$\xi_n(B(0,r))$ is ``almost round'' and has bounded degree. After this
a ``uniform'' non-normality argument is used to get that $\xi_{n+m}(B(0,r))$
eventually covers some neighbourhood $U$ of the critical points for some
$m \leq \ti{N}$, where $\ti{N}$ is ``uniform'' meaning that it only
depends on the large scale. We then show that the portion of parameters $a$ for which $\xi_{n+m}(a)$ enters
$U$ correspond to a certain fraction of the disk $B(0,r)$ for any $r > 0$ suffiently small. This will imply that the Lebeguse density of $(\de,k)$-Misiurewicz maps is stricly less than $1$ at $a=0$.

Although $\xi_n'$ has not necessarily bounded distortion on the whole
disk $B(0,r)$, one can show that we have (strong) bounded distortion of $\xi_n'$ inside so called dyadic
disks $D_0=B(a_0,r_0) \subset B(0,r)$, where $r_0/|a_0| = k <
1$ for some $k$ only depending on the family $R_a$. The sets
$\xi_n(D_0)$ will also grow to the large scale. However, to use the
``uniform'' non-normality argument we must find a point in the Julia set $\JJ(R)$
``well inside'' the set $\xi_n(D_0)$. Although one can prove
that this is the case, we choose to work directly
with the whole ball $B(0,r)$ instead of dyadic disks. (It is a matter
of taste which method one prefers. Strong argument distortion
estimates will be needed if one works with $B(0,r)$ but the proofs are
not significantly easier if one only wants the ``usual'' distortion
estimates of absolute values.)

On the other hand, we have a distortion estimate (of the absolute values) in an annular domain $A=A(r_1,r) = \{ a : r_1 < |a| < r \}$ of parameters for $\xi_n'$. These parameters will be
mapped by $\xi_n$ into an annular domain of the type
$A(\de'',\de',w)=\{ z : \de'' < |z-w| < \de' \}$.
In this case, for $a,b \in A$, the distortion estimate will be of the form
\[
\biggl| \frac{\xi_n'(a)}{\xi_n'(b)} \biggr| \leq C.
\]
The real positive numbers $\de',\de''$ where $\de' \gg \de''$ only
depend on the unperturbed function $R_0=R$.

\section{Expansion}

A fundament in getting distortion estimates is to have exponential growth of $|\xi_n'(a)|$.  The strategy we employ is similar to the one used in
\cite{MA}, Chapter 4, but the major difference here is that the postcritical
set is not (necessarily) finite. (In \cite{MA} the critical points
for the unperturbed function all land on repelling periodic orbits.) We will use a theorem by R. Ma\~n\'e to get expansion on the postcritical set for a Misiurewicz map. This expansion will then be transferred to the parameter space by Proposition \ref{vava}.


Recall that a compact set $\La$, which is invariant under $f$, is {\em
hyperbolic} if there are constants $C > 0$ and $\la > 1$ such that for any $z
\in \La$ and any $n \geq 1$,
\[
|(f^n)'(z)| \geq C \la^n.
\]

The main result which we will use by Ma\~n\'e (see \cite{RM}) is the following.

\begin{Thm}[Ma\~n\'e's Theorem I] \label{mane}
Let $f:\hat{\C} \mapsto \hat{\C}$ be a rational map and $\La \subset J(f)$ a
compact invariant set not containing critical points or parabolic points. Then either $\La$ is a hyperbolic set or $\La \cap \omega(c) \neq \emptyset$
for some recurrent critical point $c$ of $f$.
\end{Thm}
\begin{Thm}[Ma\~n\'e's Theorem II] \label{mane2}
If $x \in J(f)$ is not a parabolic periodic point and does not intersect $\om(c)$ for some recurrent critical point $c$, then for every $\vep > 0$, there is a neighborhood $U$ of $x$ such that

\begin{itemize}
\item For all $n \geq 0$, every connected component of $f^{-n}(U)$ has diameter $\leq \vep$.

\item There exists $N > 0$ such that for all $n \geq 0$ and every connected component $V$ of $f^{-n}(U)$, the degree of $f^n |_V$ is $\leq N$.

\item For all $\vep_1 > 0$ there exists $n_0 > 0$, such that every connected component of $f^{-n}(U)$, with $n \geq n_0$, has diameter $\leq \vep_1$.

\end{itemize}
\end{Thm}
An alternative proof of Ma\~n\'e's Theorem can also be found by L. Tan and M. Shishikura in \cite{ST}.

A corollary of Ma\~n\'e's Theorem II is that a Misiurewicz map cannot have any Siegel disks, Herman rings or Cremer points (see \cite{RM} or \cite{ST}).
\comm{
Indeed, take a point $x$ on the boundary of a Siegel disk or Herman ring, call it $B$, and let $U$ be a neighborhood of $x$ as in Ma\~n\'e's Theorem II. Then there is some $y \in B \cap U$, with $\dist(y,\partial B) \geq r$, for some $r > 0$. Since a Misiurewicz map has no recurrent critical points, we get that every connected component $V$ of $f^{-n}(U)$ has diameter $d_n \raw 0$ as $n \raw \infty$. Now let $V$ be such a component intersecting the boundary of $B$ and containing $f^{-n}(y)$. Inside $B$, $f$ or an iterate of $f$ is conjugate to a rotation of angle $\th \notin \Q$. Hence there is a well defined inverse $f^{-1}$ (or $f^{-p}$) inside $B$. This implies that for every $\vep_1$ there is an $n_0 > 0$ such that we can find a $n \geq n_0$ such that $|f^{-n}(y)-y| \leq \vep_1$. But then since $V=f^{-n}(U)$ contains both $f^{-n}(x) \in \partial B$ and $f^{-n}(y) \in B$ we have that the diameter of $V$ must be greater than, say, $r/2$, if $\vep_1 < r/4$. The claim follows.
}
In particular, a Misiurewicz map has no indifferent cycles.

\subsection{Expansion near the postcritical set.}

By Ma\~n\'e's Theorem, the Misiurewicz
condition gives rise to expansion of the derivative in a (closed)
neighborhood of the postcritical set. More precisely, the
postcritical set $P^k(R)$ for a Misiurewicz map $R(z)$ is hyperbolic
for some (smallest) $k > 0$. Put
$\La=P^k(R)$ for this $k$. Since $\La$ is hyperbolic there exists a
holomorphic motion $h: \La \times B(0,r) \raw
\hat{\C}$, such that for each fixed $a \in B(0,r)$, $h(z,a)=h_a: \La \raw \La_a$ is
quasiconformal and for fixed $z \in \La$ the map $h(z,a)$ is
holomorphic in $a \in B(0,r)$ (see \cite{MeSt}, Theorem III.1.6). Moreover,
\begin{equation} \label{conj}
h_a \circ R_0(z) = R_a \circ h_a(z), \text{   for all $z \in \La$}.
\end{equation}

Since $\La$ is hyperbolic there is some $N$ such that $|(R_0^N)'(z)|\geq \la_0$ for some $\la_0 > 1$ for all $z \in \La$. Now take a neighborhood $\NN$ of $\La$ such that we have expansion of $|(R_a^N)'(z)|$ for all $z \in \NN$. Thus, for some $C_1 > 0$ and $\la_1 > 1$,
\begin{equation}\label{expan}
|(R_a^N)'(z)| \geq C_1 \la_1^j,
\end{equation}
whenever $R^k(z) \in \NN$ for $k = 0,1,\dots,j$, for all $a \in B(0,r)$.
Assume moreover that $\NN$ be so that $U_{10\de} \cap \NN = \emptyset$ and that $\NN$ is closed. We get immediately the following lemma.
\begin{Lem} \label{expinside}
There exists some $\la > 1$ and $r > 0$, such that whenever $R_a^k(z) \in \NN$ for $k = 0,1,\dots,j$ and $a \in B(0,r)$, we have
\begin{equation}
|(R_a^j)'(z)| \geq C \la^j.
\end{equation}
\end{Lem}
Take some $\de' > 0$ such that $\{ z : \dist(z,\La) \leq 11\de' \} \subset \NN$. Moreover, choose $r > 0$ small enough so that
$$
\{ z : \dist(z,\La_a) \leq 10\de' \} \subset \NN, \text{ for all $a \in B(0,r)$}.
$$
Hence $\La_a$ is well inside $\NN$ for all $a \in B(0,r)$.
Further, there will be more conditions on $\de'$ in Section
\ref{distortion} (so that we might have to diminish $\de'$). Define
\[
P_{\de'} = \{ z : \dist(z,\La) < \de' \}.
\]

\subsection{The transversality condition}
For the main construction to work, we need a certain transversality
condition, meaning that the critical values must not follow {\em the
  holomorphic motion} of the critical values. Recall that every
critical point $c_j$ of $R_0$ eventually maps onto $\La$,
i.e. $R_0^{k_j}(c_j)=v_j \in \La$, for some (smallest) $k_j > 0$. These
$v_j(a)$ move holomorphically in $a$. We want to compare these
functions with the holomorphic motion of the starting critical value $v_j=v_j(0)$. Put
\begin{equation} \label{transv}
x_j(a) = v_j(a)-h_a(v_j).
\end{equation}

If all $x_j(a) \equiv 0$ then it means that every critical point
$c_j(a)$ in the Julia set is mapped onto the hyperbolic set $\La_a$ by
(\ref{conj}). Moreover, since
all the functions $\xi_{n,j}(a)$ would form a normal family, we get
that every $R_a$ would be a Misiurewicz map, quasiconformally
conjugate to $R_0$, by Theorem 4.2 in \cite{McM-book}.

Before stating the next lemma, we refer to \cite{McM-book} for the
definition of a {\em Latt\'es map}, which is a type of postcritically
finite rational map for which the Julia set is equal to the whole
Riemann sphere. These type of maps were first introduced by
S. Latt\'es \cite{SL}. We also need the two following results; the first is Teichm\"uller's module theorem (see \cite{LV}):

\begin{Thm}
Let $G$ be an annular domain which separates $0$ and $z_1$ from $z_2$
and $\infty$. Then
\[
\mod (G) \leq \log \frac{|z_1|+|z_2|}{|z_1|} + C(\frac{|z_2|}{|z_1|}),
\]
where the function $C$ is bounded by $2 \log 4$.
\end{Thm}
The second is Lemma 2.2 from \cite{JJ}:
\begin{Lem}
Let $\DD \subset \C$ be a simply connected domain and let $F:\CC \raw \D=\{ |z|> 1\}$, $F(\partial \DD) \subset \partial \D$, be $p$-valent (i.e. degree $p$). Then if $\rho$ denotes the hyperbolic metric,
\begin{align}
\{ w \in \D : \rho_{\D}(F(z_0),w) \leq C \} &\subset F( \{ z\in \DD: \rho_{\DD}(z,z_0) \leq 1 \}) \nonumber \\
&\subset \{ w \in \D : \rho_{\D}(F(z_0,w)) \leq 1 \}, \nonumber
\end{align}
where $C$ depends on $p$.
\end{Lem}



\begin{Lem} \label{pullback}
Assume that $B(0,r)$ is a good family. If there exists some $r > 0$
such that $x_j(a)=0$ for all $j$ and all $a \in B(0,r)$, then $R_0$ and $R_a$ are conformally conjugate for every $a \in B(0,r)$.

\end{Lem}
\begin{proof}
We now that any $R_a$ and $R_0$ are quasi conformally conjugate for
any $a \in B(0,r)$. We have to prove that the conjugation is in fact conformal.

If $J(R_0) \neq \hat{\C}$ then \cite{FP1} Theorem A by F. Przytycki implies that
$J(R_0)$ has Lebesgue measure zero (see also \cite{JJ} by Carleson, Jones and
Yoccoz). We now use a standard argument (see \cite{McM-book}, proof of Theorem
4.9) as follows. Since $B(0,r)$ is a good family, for any $a \in
B(0,r)$ we have that $R_a$ and $R_0$
are conformally conjugate on the Fatou set. So there exists some
conformal map $\varphi: \FF(R_0) \raw \FF(R_a)$ such that
\[
\varphi \circ R_0 = R_a \circ \varphi(z),
\]
for all $z \in \FF(R_0)$. Hence $\varphi$ is a holomorphic motion of
$\FF(R_0)$. By the $\la$-Lemma by \cite{MSS}, this motion extends to a holomorphic motion $\varphi_1$ on the
closure of $\FF(R_0)$ which must be equal to $\hat{\C}$ since
$\JJ(R_0)$ has Lebesgue measure zero. Hence $\varphi_1$ is a
quasi-conformal conjugacy
on $\hat{\C}$ which is conformal almost
everywhere. It follows that $\varphi_1$ is conformal.

If $J(R_0)=\hat{\C}$ then again we follow an argument in
\cite{McM-book} (cf. also \cite{SvS} Theorem 3.3). The
quasiconformal conjugacy induces a dilatation $\mu(z)$.
Let $E$ be the support of $\mu$ and assume that $E$ has
positive Lebesgue measure. Take a density point $z_0
\in E$ such that $z_0$ is a Lebesuge point for $\mu(z)$. Let $\th(z)$
be the angle of the invariant line field induced by $\mu(z)$.
Take a limit point $x$ of $R_0^n(z_0)$, so that $R_0^{n_k}(z_0) \raw x$.
By Ma\~n\'e's Theorem, the degree of $R_0^{n_k}$
restricted to the component $W_k$ of $R_0^{-n_k}(B(x,\eta))$
containing $z_0$ is uniformly bounded. Moreover, Ma\~n\'e's Theorem
implies that $diam(W_k) \raw 0$ as $k \raw \infty$. Since $z_0$ is a
point of density for $E$, we have
\[
\lim\limits_{k \raw \infty} \frac{m( E \cap W_k)}{m(E)} \raw 1.
\]
From Lemma 2.2 \cite{JJ} it follows that for every $C$ there is some constant $q$ only depending on $N$ and $C$ such that
\[
f^{n_k} (\{ w \in W_k : \rho_{W_k}(z,z_0)
\leq C \}) \supset \{ z \in B(x,\eta) :
\rho_{B(x,\eta)}(z,f^{n_k}(z_0)) \leq q \}.
\]
Put $W_k'=\{ w \in W_k : \rho_{W_k}(z,z_0)
\leq C \}$. Choosing $C$ sufficietly small we can ensure that
\[
\mod (W_k \sm W_k') \geq 100.
\]
Now, from Teichm\"uller's modulus theorem it follows that if $G_k = W_k \sm W_k'$ has sufficiently large modulus, $\mod G \geq 100$ will do, then there is a ball $B_k$
with boundary $C_k$, such that $C_k \subset G_k$.
For every $k$, we have
\[
R_0^{n_k}(B_k) \supset \{ z \in B(x,\eta) :
\rho_{B(x,\eta)}(z,R_0^{n_k}(z_0)) \leq q \}.
\]

Let $A_k:B_k \raw
D$ be a linear normalisation of $B_k$, where $D$ is the unit disk.
Then $f_k = R_0^{n_k} \circ A_k^{-1} : D
\raw B(x,\eta)$ is a normal family and $$f_k(D) \supset \{ z \in B(x,\eta) :
\rho_{B(x,\eta)}(z,f^{n_k}(z_0)) \leq q \}.$$ Hence there is a subsequence
$k'$ such that $f_{k'}$ converges uniformly to a non-constant limit
function $f$. Since $z_0$ is a Lebesgue point for $\th$, the family of
line fields $(A_k)_*(\th)$ tends to a constant line field $\th'$ in
$D$. Hence $f_*(\th')$ is a holomorphic line field which coincides
with $\th$.

Lemma 3.16 in \cite{McM-book} now implies that $R_0$ has to be a Latt\'es map or $E$
has measure zero. In other words, $R_a$ and $R_0$ are conformally
conjugate if they are not Latt\'es maps. If they are Latt\'es maps,
then Thurston's Theorem implies that $R_0=R_a$, since $R_0$ and $R_a$ are quasi-conformally conjugate. The lemma follows.

\end{proof}

So if all $x_j(a)$ are identically equal to zero then in fact the
maps in $B(0,r)$ are all conformally conjugate Misiurewicz maps, which
is impossible since the family $f_a$, $a \in B(0,r)$ is assumed to be normalised. Hence we have the following important transversality criteria:

{\bf Transversality criteria}:
{\em There is at least one critical value $v_j(0) \in \La$, where $\La$ is
  a holomorphically moving hyperbolic set, for which
  $x_j(a)=v_j(a)-h_a(v_j(0))$ is not identically equal to zero.}


For this reason, by $\xi_n(a)$ we mean
$\xi_{n,j}(a)$ for this particular $j$ for which the transversality condition holds further on, unless otherwise stated.

\section{Distortion Lemmas} \label{distortion}
This section is devoted to the distortion results, which are used to get control of $\xi_n(B(0,r))$ up to the large scale. One of the main result in this section is Proposition \ref{vava}, which shows a close relation between the space derivative and parameter derivative. The idea of comparing these two quantities was first introduced by M. Benedicks and L. Carleson in \cite{BC2} (see also \cite{BC1}). We assume throughout this section that the transversality criteria holds.

By the transversality criteria there is some critical value $v_j(a)$ for which $x_j(a)$ is not identically equal to zero. Hence putting $x(a)=x_j(a)$ we have
\begin{equation} \label{xK1}
x(a) = K_1 a^k + \ldots,
\end{equation}
for some $K_1 \neq 0$. Define $\mu_n(a)=\mu_{n,j}(a)=h_a(R_0^n(v_j))$.
Then in particular $x(a) = \xi_0(a) - \mu_0(a)$.

\comm{
Given $\de' > 0$, in the following, $\de'' > 0$ shall be thought of as a
small fixed number, satisfying $|\log \de''| \ll
|\log |x(r)||$ (this holds if the  perturbation $r > 0$ is small
enough), and $\de'' < \de'/M_0^{10}$, where $M_0=\sup |R_a'(z)|$, where the supremum is taken over all $z \in \hat{\C}$ and $a \in B(0,r)$ (there is a lot of freedom of how to choose $\de''$, roughly speaking it should be significantly smaller than $\de'$).
}

The hyperbolic set $\La$ and its neighbourhood $\NN$ will be the backbone in the expansion and distortion estimates. Let us first state an elementary property, saying that two points close two each other inside $\NN$ repel each other uniformly up to some large scale. By the definition of $\NN$ there exists constants $N > 0$ and $\la > 1$ such that $|(R_a^N)'(z)| \geq \la$ for all $z \in \NN$ and $a \in B(0,r)$ for some $\la > 1$.

Hence to every $z \in \NN$ there is some radius $r(z) > 0$ such that
\begin{equation} \label{hyper}
|R_a^N(z)-R_a^N(w)| \geq \la |z-w|,
\end{equation}
for all $w \in \NN$ satisfying $|z-w| \leq r(z)$ (possibly diminishing $\la > 1$ slightly). Since $\NN$ is compact and $r(z)$ is continuous there is a constant $\ti{r} > 0$ such that (\ref{hyper}) holds for all $z,w \in \NN$ provided $|z-w|\leq \ti{r}$. For simplicity assume that $N=1$.

The following lemma, which will be needed in the subsequent lemma, is variant of Lemma 15.3 in \cite{WR} (see also \cite{MA}, Lemma 2.1).
\begin{Lem} \label{prod-dist}
Let $u_n \in \C$ be complex numbers for $1 \leq n \leq N$. Then
\begin{equation}
\biggl| \prod_{n=1}^N (1+u_n)-1 \biggr| \leq
\exp \biggl(\sum_{n=1}^N |u_n| \biggr) - 1.
\label{ineq-1}
\end{equation}
\end{Lem}
In the following, by $(R_a^n)'(\mu_0(a))$ or $(R_a^n)'(v(a))$ we mean $(R_a^n)'(z)$ evaluated at $z=\mu_0(a)$ or $z=v(a)$ respectively.
\begin{Lem}[Main Distortion Lemma] \label{distortion1}
For every $\vep > 0$, there are arbitrarily small constants $\de' > 0$ and $r > 0$ such that the following holds. Let $a,b \in B(0,r)$ and suppose that $|\xi_k(t) - \mu_k(t)| \leq \de'$, for $t=a,b$ and all $k \leq n$. Then
\begin{equation} \label{distest}
\biggl| \frac{(R^n)'(v(a),a)}{(R^n)'(v(b),b)} - 1 \biggr| < \vep.
\end{equation}
The same statement holds if one replaces $v(s)=\xi_0(s)$, $s=a,b$, by $\mu_0(t)$, $t=a,b$ in (\ref{distest}).
\end{Lem}

\begin{proof}
The proof goes in two steps. Let us first show that
\begin{equation}
\biggl| \frac{(R_t^n)'(\mu_0(t))}{(R_t^n)'(\xi_0(t))} - 1 \biggr| \leq \vep_1,
\label{close-to-1}
\end{equation}
where $\vep_1=\vep(\de')$ is close to $0$. We have
\begin{align}
\sum_{j=0}^{n-1} \biggl| \frac{R_t'(\mu_j(t)) - R_t'(\xi_j(t))}{R_t'(\xi_j(t))} \biggr| &\leq C_{\de} \sum_{j=0}^{n-1} |R_t'(\mu_j(t)) - R_t'(\xi_j(t))| \nonumber \\
&\leq C_{\de}\max|R_a'(z)| \sum_{j=0}^{n-1} |\mu_j(t)-\xi_j(t)| \nonumber \\
&\leq C \sum_{j=0}^{n-1} \la^{j-n}|\mu_n(t)-\xi_n(t)| \leq C(\de'), \nonumber
\end{align}
where we used equation (\ref{hyper}). By Lemma \ref{prod-dist}, (\ref{close-to-1}) holds if $\de'$ is small enough.
Secondly, we show that
\[
\biggl| \frac{(R_t^n)'(\mu_0(t))}{(R_s^n)'(\mu_0(s))} - 1\biggr|\leq \vep_2,
\]
where $\vep_2=\vep_2(\de')$ is close to $0$. Put
$\la_{t,j}=R_t'(\mu_j(t))$. Since $\la_{t,j}$ are all analytic in $t$ we have
$\la_{t,j}=\la_{0,j}(1+c_j t^l + \ldots)$. Moreover, since $n \leq -C \log|x(t)| = -C \log |t|$
\[
\frac{(R_t^n)'(\mu_0(t))}{(R_s^n)'(\mu_0(s))}  = \prod \frac{\la_{t,j}}{\la_{s,j}} = \prod_{j=0}^{n-1}
\frac{\la_{0,j}(1+c_j t^l + \ldots)}{\la_{0,j}(1+c_j s^l + \ldots)} =
\frac{1+cn t^l + \ldots}{1+cn s^l + \ldots}.
\]
Both the last numerator and denominator in the above equation can be
estimated by $1 + \OO((\log |t|) |t|^l)$ and $1 + \OO((\log |s|)
|s|^l)$, which both can be made arbitrarily close to $1$ if $r > 0$ is
small enough. From this the lemma follows.
\end{proof}

\begin{Lem} \label{xinprim}
Let $\vep > 0$. If $\de' > 0$ is sufficiently small, then for every $0 < \de'' < \de'$ there exist $r > 0$ such that the following holds. Let $a \in B(0,r)$ and assume that $|\xi_k(a)-\mu_k(a)| \leq \de'$, for all $k \leq n$ and $|\xi_n(a)-\mu_n(a)| \geq \de''$. Then
\[
\biggl| \frac{\xi_n'(a)}{(R_a^n)'(\mu_0(a))x'(a)} - 1 \biggr| \leq \vep.
\]
\end{Lem}
\begin{proof}
First we note that by Lemma \ref{distortion1} we have
\[
\xi_n(a) = x(a) (R_a^n)'(\mu_0(a)) + \mu_n(a) + E_n(a),
\]
where, for instance $|E_n(a)| \leq |\xi_n(a) - \mu_n(a)|/1000$
independently of $n$ and $a$ if $\de'$ is small enough. Put $R_a'(\mu_j(a)) = \la_{a,j}$. Differentiating with respect to $a$ we get
\begin{equation} \label{xip}
\xi_n'(a) = \prod_{j=0}^{n-1} \la_{a,j} \biggr( x'(a) + x(a) \sum_{j=0}^{n-1} \frac{\la_{a,j}'}{\la_{a,j}} + \frac{\mu_n'(a) + E_n'(a)}{\prod_{j=0}^{n-1} \la_{a,j}}  \biggl).
\end{equation}

We claim that only the $x'(a)$ is dominant in (\ref{xip}) if $n$ is large so that $\de'' \leq |\xi_n(a)-\mu_n(a)| \leq \de'$. This means that, by Lemma \ref{distortion1},
\[
(1-\vep_1)\de'' \leq |x(a)| \prod_{j=0}^{n-1} |\la_{a,j}| \leq (1+\vep_1) \de' < 1,
\]
where $\vep_1 > 0$ is arbitrarily small provided $r > 0$ is small enough.
Since $\prod_{j=0}^{n-1} |\la_{a,j}| \geq \la^n$, for some $\la > 1$, taking logarithms and rearranging we get
\begin{equation}
(1-\vep_1) \sum_{j=0}^{n-1} \log |\la_{a,j}| \leq -\log |x(a)| \leq (1+\vep_1) \sum_{j=0}^{n-1} \log |\la_{a,j}|, \label{nsize}
\end{equation}
if $|\log \de''|  \ll |\log |x(a)||$, which is true if the perturbation $r > 0$ is chosen sufficiently small compared to $\de''$.
Since $|\la_{a,j}| \geq \la > 1$, this means that
\[
|x(a)| \sum_{j=0}^{n-1} \frac{|\la_{a,j}'|}{|\la_{a,j}|} \leq |x(a)| n C \leq -C |x(a)| \log|x(a)|.
\]
Finally $-|x(a)| \log|x(a)| / |x'(a)| \raw 0$ as $a \raw 0$.

Now, $|E_n(a)|$ is uniformly bounded in $B(0,r)$. Therefore,
$|E_n'(a)|$ is also uniformly bounded on compact subsets of
$B(0,r)$ by Cauchy's Formula. By diminishing $r > 0$ slightly we can assume that
both $|E_n(a)|$ and $|E_n'(a)|$ are uniformly bounded on $B(0,r)$. Hence, the last two
terms in (\ref{xip}) tend to zero as $n \raw \infty$, since also
$|\mu_n'(a)|$ is uniformly bounded. We have proved that
\[
\biggl| \xi_n'(a) - x'(a) \prod_{j=0}^{n-1} \la_{a,j} \biggr| \leq \vep |\xi_n'(a)|,
\]
if $|\xi_n(a)-\mu_n(a)| \leq \de'$ and $n \geq N$ for some $N$. Choose the perturbation $r$ sufficiently small so that this $N$ is at most the number $n$ in (\ref{nsize}). Since $\la_{a,j} = R_a'(\mu_j(a))$, the proof is finished.
\end{proof}

Combining Lemma \ref{xinprim} and Lemma \ref{distortion1} we immediately get the following important result.
\begin{Prop} \label{vava}
Let $\vep > 0$. If $\de' > 0$ is small enough and $0 < \de'' < \de'$, there is an $r > 0$ such that the following holds. Take any $a \in B(0,r)$ and assume that $|\xi_k(a)-\mu_k(a)| \leq \de'$, for all $k \leq n$ and $|\xi_n(a) - \mu_n(a)| \geq \de''$. Then
\begin{equation}
\biggl| \frac{\xi_n'(a)}{(R_a^n)'(v(a)) x'(a)} - 1 \biggr| \leq \vep.
\end{equation}
\end{Prop}

\subsection{Distortion in an annular domain.} \label{annulus}
Since we may have $x'(0)=0$, we have to stay away from $a=0$ in order
to hope for distortion estimates on $\xi_n'$ on the ball $B(0,r)$. We therefore consider an annular domain $A=A(r_1,r)=\{a : r_1 < |a| < r \}$ in the parameter disk $B(0,r)$, for fixed $0 < r_1 < r$. For $a,b \in A$ we see that
\[
\frac{1}{C} \biggl( \frac{r_1}{r} \biggr)^{k-1} \leq \frac{|x'(a)|}{|x'(b)|} \leq C \biggl( \frac{r}{r_1} \biggr)^{k-1},
\]
for some constant $C \geq 1$.
From Lemma \ref{distortion1} and Proposition \ref{vava} it follows that
\begin{equation}\label{globalxi}
\frac{1}{C'} \biggl( \frac{r_1}{r} \biggr)^{k-1} \leq \frac{|\xi_n'(a)|}{|\xi_n'(b)|} \leq C' \biggl( \frac{r}{r_1} \biggr)^{k-1},
\end{equation}
for some constant $C' \geq 1$, for all $a,b \in A$ as long as $\de'' \leq |\xi_n(a)-\mu_n(a)| \leq \de'$, for all $a \in A$.

\begin{Lem} \label{annulus-l}
Let $\vep > 0$ . Then if the numbers $\de' > 0$ and $\de''/\de'$ are sufficiently
small where $0 < \de'' < \de'$, there exists an $r > 0$ such that the
following holds for any ball $B=B(0,r_2) \subset B(0,r)$.

Assume that $n$ is maximal for which $|\xi_n(a)-\mu_n(a)| \leq \de'$ for all $a \in B$. Let $r_1 < r_2$ be minimal so that $|\xi_n(a)-\mu_n(a)| \geq \de''$ for all $a \in A=A(r_1,r_2)$. Then $r_1/r_2 < 1/10$ and there is some $\de_1' \leq \de' $ such that $\xi_n(A) \subset A(\de''-\vep,\de_1'+\vep,\mu_n(0))$ and $A(\de''+\vep,\de_1'-\vep,\mu_n(0)) \subset \xi_n(A)$.

Moreover, $\xi_n$ is at most $k$-to-$1$ on $B$.
\end{Lem}
\begin{proof}
First note that a small parameter circle $\ga_r = \{a \in B(0,r): |a| = r\}$ is mapped under $x(a)$ onto a curve that encircles $0$ $k$ times, and such that $x(\ga_r)$ is arbitrarily close to a circle of radius $K_1r^k$. We have that $|\xi_n(a)-\mu_n(a)| \gg |\mu_n(a)-\mu_n(0)|$, if $|\xi_n(a)-\mu_n(a)| \geq \de''$. By Lemma \ref{distortion1}, for every $\vep_1 > 0$ it is possible to choose $r > 0$ and $\de' > 0$ such that
\begin{equation} \label{circles}
| \xi_n(a)-\mu_n(a) - (R_a^n)'(v(a))x(a) | \leq \vep_1 |\xi_n(a)-\mu_n(a)|,
\end{equation}
for all $a \in B$. Moreover, since (\ref{circles}) holds for all $a \in B$, Lemma \ref{distortion1} implies that $$\de'/\de'' \leq C|x(a_2)|/|x(a_1)|,$$ for some constant $C \geq 1$ (arbitrarily close to $1$) for all $a_1,a_2 \in B$, where $|a_1|=r_1$ and $|a_2|=r_2$. Hence we can easily choose $\de'' > 0$ so that $\de'/\de''$ is so large so that $r_1 /r_2 \leq 1/10$. From this the first part of the lemma follows.

Moreover, as the parameter $a$ orbits around the circle $\ga_r$ once,
Lemma \ref{distortion1} and (\ref{circles}) implies that $\xi_n(a)-\mu_n(a)$ orbits around $0$ $k$ times (note that here strong argument distortion estimates is needed). Since $|\mu_n(a)-\mu_n(0)| \ll |\xi_n(a)-\mu_n(a)|$, this means that $\xi_n(a)$ orbits around $\mu_n(0)$ $k$ times close to a circle of radius  $|\xi_n(a)-\mu_n(0)|$ centered at $\mu_n(0)$. By the Argument Principle, the map $\xi_n$ is at most $k$-to-$1$ on $B$.
\end{proof}

Hence an annulus $A(r_1,r_2) \subset B(0,r)$ is mapped onto a slightly distorted annulus $\xi_n(A)$ inside $A(\de'',\de', \mu_n(0))$, and the set $\xi_n(B(0,r_2))$ is an almost round ball.
\comm{
For any $a \in B(0,r)$, we say that $\xi_n(a)$, or simply the
parameter $a$ itself, has {\em escaped} at time $n$ if $|\xi_n(a)-\mu_n(a)| \geq
\de''$. Thus for any escape situation, we have that $\xi_n'(a) \sim
x'(a) (R_a^n)'(v(a))$. Apparently, the annulus $A$ above has escaped (at
time $n$).
}
\comm{
We conclude that there is a minimal $r_1=r_1(\de'')$ such that
$A(r_1,r)$ grows almost uniformly up to the large scale so that
$\xi_n(A) \subset A(\de'',\de',\mu_n(0))$, where $n$ is chosen maximal so that
$|\xi_n(a) - \mu_n(a)| \leq \de'$, for all $a \in
B(0,r)$.
}


\section{The free period} \label{free}
The main object of this section is to show that once the set
$\xi_n(B)$, for some ball $B=B(0,r_2) \subset B(0,r)$, has reached
diameter $\de' > 0$ then $\xi_{n+m}(B)$ will cover $\oli{U}_{3\de/4}$ within a finite number of iterates, i.e. $m \leq \ti{N}$ for some $\ti{N}$ only depending on $\de'$. These last $m$ iterates are referred to as the {\em free period}. We begin this section with the following elementary and important lemma.


\begin{Lem} \label{finret}
For every $d > 0$ there is an $r > 0$ such that the following holds. Let $D$ be a set which contains a disk of radius $d$ centered at the Julia set of $R$. Then there is some integer $\ti{N}$ only depending on $R$ and $d$ such that
\[
\inf \{ m \in \N : R^m(D) \supset
\oli{U}_{\de/2} ) \} \leq \ti{N}.
\]
\end{Lem}
\begin{proof}
First cover $\JJ(R)$ with the collection of open
disks $D_z$ of diameter $d$ centered at any point $z$ of $\JJ(R)$.
Since $R^n$ is not normal on the Julia set, we get that for every $D_z$, there is a smallest number $n=n(z)$ such that $R_0^n(D_j) \supset \oli{U}_{\de/2}$. Note that $n(z)$ is constant in some neighbourhood of $z$ since $R^n$ is a continuous function. Since $\JJ(R)$ is compact there is some uniform $\ti{N}$ such that $n(z) \leq \ti{N}$. The lemma is proved.
\end{proof}
Since $R_a^n(D_z)$ moves continuously in $a$, there
is an $r > 0$ such that the same statement holds for $R_a$ instead of
$R_0$, if $a \in B(0,r)$. Moreover, if $d$ depends only on $\de'$, note that $\ti{N}$ depends only on $\de'$.

Let $B \subset B(0,r)$ be a ball centered at $0$ and
assume that $n$ is maximal such that $\xi_n(B)$ has diameter at most $\de'$. It follows from Lemma \ref{annulus-l} that the set $\xi_n(B)$ contains a ball $D_1$ of diameter $d=\de'/(2M_0)$, centered at $\xi_n(0) \in J(R_0)$,
where $M_0=\max |R_a'(z)|$ for $z \in \hat{\C}$, $a \in B(0,r)$.
Lemma \ref{finret} now gives, setting $D=\xi_n(B)$, that $R_0^m(D) \supset \oli{U}_{\de/2}$, for some minimal $m \leq \ti{N}$.

We now estimate the measure of those $z \in D$ that ge mapped into $U$
under $R_0^j$, for some $j \leq m$. The following lemma may seem
evident; the only point being that the constant $C$ ``a priori'' does
not depend on $D$.

\comm{
Let us keep in mind the collection $\DD$ of disks $D_z$ in Lemma \ref{finret} of diameter $d$ centered at points $z$ in the Julia set of $R$.
\begin{Lem}
There is some uniform $C > 0$ only depending on $d$, $f$ such that for any $D_w \in \DD$, we have
\[
\mu(\{ z \in D_w: R^j(z) \in U, \text{ for some $0 \leq j \leq m$} \})
\geq C \mu (D_w).
\]
\end{Lem}
\begin{proof}
By Lemma \ref{finret}, $m \leq \ti{N}$, where $\ti{N}$ only depends on $d$. Moreover for each $D_z \in \DD$ the measure $\mu(W_z)$ of the preimage $W_z=R^{-m}(U) \cap D_z$ is continuous in $z$ and positive for all $z \in \JJ(R)$. Hence since $\JJ(R)$ is compact, there is some $C_0 > 0$ such that $\mu(W_z) \geq C_0$ for all $z \in \JJ(R)$.
\end{proof}
}
\begin{Lem} \label{unimeas}
Let $f: D \raw \hat{\C}$ be a rational map,
$D$ is an open set, $U$ is a neighbourhood of $Crit(f)$ and $U \cap D = \emptyset$. Assume that $f^m(D) \supset \oli{U}$ for some $m > 0$. Then there exists a constant $C > 0$ only depending on $U$, $m$ and $f$ such that
\[
\mu(\{ z \in D: f^j(z) \in U, \text{ for some $0 \leq j \leq m$} \})
\geq C \mu (D).
\]
\end{Lem}
\begin{proof}
Put
$E = \{ z \in D: f^m(z) \in U \}.$
We have
\[
\mu(U) \leq \int_{E}
|(f^m)'(z)|^2 d \mu(z) \leq C_1 \mu(E),
\]
where $C_1$ only depends on $f$, $U$ and $m$. Now let $g(w) = \{z \in D: f^m(z)=w \}$. We get
\[
\mu(D \sm E) = \int_{\hat{\C} \sm U} \sum_{z\in g(w)} \frac{1}{|(f^m)'(z)|^2} d \mu(w) \leq C_2 \mu(\hat{\C} \sm U),
\]
where $C_2$ only depends on $f$, $U$ and $m$.
Since $\mu(U) \geq C_3 \mu(\hat{\C} \sm U)$ for some $C_3 > 0$ only
depending on $U$, we get
\begin{align}
\mu(E) &\geq C_1 \mu(U) \geq C_1 C_3 \mu(\hat{\C} \sm U) \geq \frac{C_1
  C_3}{C_2} \mu(D \sm E) \geq C \mu(D),
\end{align}
for some constant $C > 0$ only depending on $m$ and $U$. Since $$E \subset \{ z \in D: f^j(z) \in U, \text{ for some $0 \leq j \leq m$} \}$$ the lemma follows.
\end{proof}

Let $A(z)$ be the set of preimages in $B(0,r)$ of $z \in D$ under $\xi_n$. We arrive at the following Proposition, which we state in a more general context.
\begin{Prop} \label{meas}
Let $f_a$, $a \in B=B(0,\vep)$ be an analytic family of rational maps for some $\vep > 0$. Assume that the transversality criteria is satisfied for some critical value $v_j(a)=v(a)$. Then there exists some $\de' > 0$ and $0 < r_1 \leq \vep$ only depending on $f_0$ such that for any $0 < r < r_1$, there is some maximal $n < \infty$ such that
\[
diam(\xi_n(B)) \leq \de',
\]
and such that $\xi_n(B)$ contains a ball of diameter $\de'/(2M_0)$ centered at $\JJ(f)$. The degree of $\xi_n(B)$ is bounded by some constant $K < \infty$ only depending on the family $f_a$, $a \in B(0,r)$.

Moreover, if $U$ is an open set, and $D=\xi_n(B)$, there are constants
$C > 0$ and $\ti{N}$ such that
\begin{align}
\mu(\{ &z \in D: \xi_{n+j}(a(z)) \in U, \text{ for all $a(z) \in
  A(z)$ some $0 \leq j \leq \ti{N}$ } \})  \nonumber \\
  & \geq C \mu (D), \nonumber
\end{align}
where $C$ only depends on $f_0$ and $U$.
\end{Prop}
\begin{proof}
The first part of the proposition follows from Lemma
\ref{annulus-l}. To prove the second part, by Lemma
\ref{finret}, there is some $\ti{N} < \infty$ depending only on $\de'$ such
that $\oli{U} \subset f_0^m(D)$ for some $m \leq \ti{N}$. Let us apply
Lemma \ref{unimeas} to $f_0$ and some $U' \Subset U$. We get that
\[
\mu(\{z \in D: f_0^j(z) \in U' \text{ for some $0 \leq j \leq \ti{N}$} \})
\geq C \mu(D),
\]
for some $C$ only depending on $f_0$, $U'$ and $\de'$. Since $\de' > 0$ only
depends on $f_0$, $C$ only depends on $f_0$ and $U'$.
Since the parameter dependence can be made arbitrarily small under the
iterates from $n$ up to $n+m$ where $m \leq \ti{N}$, we get
immediately that for any open subset $U' \Subset
U$ there is some $r_1 > 0$ such that
any parameter $a \in B(0,r_1)$ for which $f_0^j(\xi_n(a)) \in U'$, has
that $\xi_{n+j}(a) \in U$. This proves the second part of the proposition.
\end{proof}
\comm{
\begin{Lem}
Assume that the diameter of $\xi_n(B)$ is at most $\de'$ for the largest possible $n$, and that $m$ is minimal so that $\xi_{n+m}(B) \cap U_{\de/10} \neq \emptyset$. Then the degree of $\xi_{n+m}$ on $B$ is bounded by some $M < \infty$, regardless of $n$.
\end{Lem}

\begin{proof}
Consider the almost round disk $D=\xi_n(B)$. The degree of $R_0^m$ on $D$ is bounded by $d^{\ti{N}}$. Moreover, for any point $z \in D$, by definition we have that $R_a^j(z)$ does not intersect $U_{\de/10}$ for $j=0,\ldots, m-1$. This means that $|(R_a^m)'(z)| \geq C $ for some $C=C(\de,\ti{N})$ and there are no critical points of $R_a^m$ inside $D$. It follows that any two points $z_1,z_2 \in D$ mapped onto the same point must be separated by some fixed constant $c$. Now, to switch from $R_a^m(D)$ to $\xi_{n+m}(B)$ we note that each pair of points $z_1,z_2 \in D$ are images under $\xi_n(a)$ for some $a \in B$, i.e. $z_1 = \xi_n(a_1), z_2 = \xi_n(a_2)$. The parameter dependence under the coming $m \leq \ti{N}$ iterates can be made arbitrarily small if the perturbation is sufficiently small (i.e if $r > 0$ is sufficiently small). Hence there is a slightly smaller constant $c/2$ such that if $z_1=\xi_n(a_1), z_2 = \xi_n(a_2)$ then $\xi_{n+m}(a_1)=\xi_{n+m}(a_2)$ only if $|z_1-z_2| \geq c/2$ or $z_1=z_2$. Since the map $\xi_n$ is at most $k$-to-$1$, the lemma follows.
\end{proof}
}

\section{Conclusion and proof of Theorem \ref{discthm}}
Assume that $R_0=R$ satisfies the Misiurewicz condition and that $B(0,r)$ is a good family around $R$. We want to prove that the set of $(\de,k)$-Misiurewicz maps in $B(0,r)$ has Lebesgue density strictly less than $1$ at $a=0$. The proof consists of showing that a specific fraction $q > 0$ of parameters in any disc $B=B(0,r_2) \subset B(0,r)$
corresponds to functions $R_a$ which have a critical point which
returns into $U_{3\de/4}$ after more than $k$ iterations. Since we know by Lemma \ref{pullback} that the transversality criteria is fulfilled for some critical point $c(a)$, we will study the iterates of this particular critical point. Put $\xi_n(a)=R_a^n(v(a))$, where $v(a)=R_a^{k_0}(c(a))$ for some smallest $k_0 > 0$ such that $R_0^{k_0}(c(0))$ does not contain any critical point in its forward orbit.

Let us assume that $N$ is largest positive integer such that $diam(\xi_{N}(B)) \leq \de'$. Moreover, choose the perturbation $r$ sufficiently
small so that $N > k$. Put $D=\xi_N(B)$.
Now we turn to the annular domain $A=A(r_1,r_2)$ as in Subsection \ref{annulus}, and the ball $B=B(0,r_2) \subset B(0,r)$.

\begin{Claim}
It is possible to choose the constant $\de'' > 0$ such that
$0 < \de'' < \de'$ and such that any $a \in B(0,r_2)$ for which $\xi_{N+j}(a) \in U_{3\de/4}$ has that $a \in A(r_1,r_2)$, where $A(r_1,r_2)$ is the
annular domain in Lemma \ref{annulus-l}.
\end{Claim}
To prove the Claim, we note that in Proposition \ref{vava} we are free to choose $\de'' > 0$ as small as desired, provided $r > 0$ is chosen small enough. It means that if $\de''$ is sufficiently small, then a parameter $a \in B(0,r)$ for which $|\xi_N(a) - \mu_N(a)| \leq \de''$ satisfies $$|\xi_{N+j}(a)-\mu_{N+j}(a)| \leq |(R_a^j)'(\xi_N(a))||\xi_N(a)-\mu_N(a)| \leq \de'$$ for all $j \leq \ti{N}$, if $|(R_a^{\ti{N}})'(z)| \leq \de'/\de''$ for all $z \in \NN$. From this the claim follows.

Moreover, recall that we have bounded distortion inside $A=A(r_1,r_2)$:
\begin{equation} \label{xidist}
\frac{1}{C} \leq \frac{1}{C'} \biggl( \frac{r_1}{r_2} \biggr)^{k-1} \leq \biggl| \frac{\xi_{N}'(a)}{\xi_{N}'(b)} \biggr| \leq C' \biggl( \frac{r_2}{r_1} \biggr)^{k-1} \leq C,
\end{equation}
for any $a,b \in A(r_1,r_2)$, where $r_1/r_2 \leq 1/10$ according to Lemma \ref{annulus-l}. Hence $\mu(A) \geq (99/100) \mu(B)$. Put
\[
E = \{ z \in D: \xi_{N+j}(a(z)) \in U_{3\de/4} \text{ for all $a(z)
  \in A(z)$ and some $j \leq \ti{N}$} \}.
\]
By Proposition \ref{meas} we have $\mu(E) \geq C \mu(D)$. Recall that
$\xi_N$ has bounded degree on $B$ by Lemma \ref{annulus-l}. By
(\ref{xidist})
\[
\mu(\{a \in B(0,r): \xi_{N}(a) \in E \}) \geq q \mu(B(0,r)),
\]
for some $0 < q < 1$, only depending on $U$, $\de'$. By the definition of $E$,
\[
\mu(\{a \in B(0,r): \xi_{n}(a) \in U_{3\de/4} \text{ for some $n > k$ }\}) \geq q \mu(B(0,r)),
\]
We now make a minor note that it may happen that for some $a \in B(0,r)$, a critical point $c_i(a) \in U'$ lies in a super-attracting cycle, where $U'$ is a component of $U$. Hence in particular, if the critical point $c(a)$ returns into this neighbourhood $U'$, the map $R_a$ may still be $(\de,k)$-Misiurewicz. However, such super-attracting cycle cannot be persistent in $B(0,r)$ since $c_i(0) \in \JJ(R_0)$. Hence the Lebesgue measure of those $a \in B(0,r)$ for which some $c_i(a) \in U$ lies in a super-attracting cycle is zero.

Since the critical points move inside $U_{\de^{10}}$ as $a \in
B(0,r)$, this will ensure that for almost all $a$ the following is true:
If a critical point $c(a)$ returns into a slightly smaller $U_{3\de/4} \subset U_{\de}$, $R_a$ cannot be $(\de,k)$-Misiurewicz. In other words,
\[
\mu(\{a \in B(0,r): \text{ $R_a$ is not $(\de,k)$-Misiurewicz} \} ) \geq q \mu(B(0,r)).
\]
Since this is true for every arbitrarily small $r > 0$, the Lebesgue
density of $(\de,k)$-Misiurewicz maps at $a=0$ is at most $1-q <
1$. The proof of Theorem \ref{discthm} is finished.

\bibliographystyle{plain}
\bibliography{ref}

\comm{
\bibliographystyle{plain}
\bibliography{bib}

}

\end{document}